\documentclass[12pt]{article}

\usepackage{amsmath,amsthm}
\usepackage{amssymb}
\usepackage{color}
\usepackage{xspace}
\usepackage[colorlinks=true,
linkcolor=green,
filecolor=brown,
citecolor=green]{hyperref}
\def\red{\textcolor{red} }

\def\s{3-$\underline{5}$-2-4-1-satisfying\xspace}
\def\m{$31$-\,$4$-$2$-avoiding\xspace}

\def\lr{LRmax\xspace}
\def\red{\textcolor{red} }

\begin{document}
\newtheorem{lemma}{Lemma}
\newtheorem*{theorem}{Theorem}
\newtheorem{prop}{Proposition}
\newtheorem{cor}{Corollary}
\begin{center}
{\Large
A Wilf Equivalence Related to Two Stack Sortable Permutations                           \\ 
}

\vspace{10mm}
DAVID CALLAN  \\
Department of Statistics  \\
\vspace*{-2mm}
University of Wisconsin-Madison  \\
\vspace*{-2mm}
Medical Science Center \\
\vspace*{-2mm}
1300 University Ave  \\
\vspace*{-2mm}
Madison, WI \ 53706-1532  \\
{\bf callan@stat.wisc.edu}  \\
\vspace{5mm}

October 10, 2005
\end{center}

\vspace{3mm}
\begin{center}
   \textbf{Abstract}
\end{center}
A permutation is so-called two stack sortable if it (i) avoids the (scattered) 
pattern \mbox{2-3-4-1}, and (ii) contains a \mbox{3-2-4-1} pattern only as part of a \mbox{3-5-2-4-1} pattern.
Here we show that the permutations on $[n]$ satisfying condition (ii) 
alone are equinumerous with the permutations on $[n]$ that avoid the 
mixed scattered/consecutive pattern \mbox{31-4-2}. The proof uses a known 
bijection from \mbox{3-2-1}-avoiding to \mbox{3-1-2}-avoiding permutations.

\vspace{10mm}

{\Large \textbf{1 \quad Introduction}  }

An instance of the pattern 31-4-2 in a permutation on $[n]$, 
considered as a list, is a sublist $abcd$ whose reduced form 
(replace smallest entry by 1, next smallest by 2 and so on) is 3142 
and in which $a$ and $b$ are consecutive entries in the permutation. 
Similarly for other patterns, the absence of a dash mandating 
adjacency in the permutation. Thus, for 
example, 35142 fails to be \mbox{3-1-4-2}-avoiding but is \mbox{31-4-2}-avoiding.
We say a permutation that  contains a 3-2-4-1 pattern only as part of a 
3-5-2-4-1 pattern is 3-$\underline{5}$-2-4-1-satisfying. See 
\cite{eigen} for a characterization of the generating function for \s 
permutations.

Our main result is
\begin{theorem}
    The number of $3$-$\underline{5}$-$2$-$4$-$1$-satisfying permutations on $[n]$ is the same as the number 
    of $31$-\,$4$-$2$-avoiding permutations on $[n]$.
\end{theorem}

Permutations avoiding one or two dashed patterns on 3 letters have been 
counted \cite{claesson}. Some generating functions and asymptotic 
results for permutations avoiding consecutive patterns, along with 
bounds for some dashed 4-letter patterns,  
are given in \cite{consecutive,asymptotic05}. 

\vspace{10mm}

{\Large \textbf{2 \quad Proof of Theorem}  }

To prove the theorem, we use a recursively-formulated 
characterization, involving left to right maxima, for  
each of the two classes of permutations in question. The left to right 
maximum (\lr ) entries of a permutation are the ``record highs'' 
when the permutation is scanned from left to right. Let 
$M=(m_{i})_{i=1}^{r}$ denote the \lr entries and  
$P=(p_{i})_{i=1}^{r}$ their respective positions. For example, the 
permutation 
$\ \red{\overset{1}{3}}\:\overset{2}{1}\:\red{\overset{3}{5}}\:\overset{4}{4}\:
\overset{5}{2}\:\red{\overset{6}{7}}\:\overset{7}{6}\ $ 
has $r=3,\ P=(1,3,6),\ M=(3,5,7)$.
Say the pair $(P,M)$ is the \emph{\lr specification} of a permutation. 
The valid \lr specifications for $[n]$ (meaning there exists a permutation 
on $[n]$ with that specification) are easily seen to be 
characterized by the conditions
\[
\begin{array}{c}
    1=p_{1}<p_{2}<\ldots<p_{r}\le n,  \\
    1\le m_{1}<m_{2}<\ldots<m_{r}=n, \\
    p_{i+1}\le m_{i}+1 \textrm{\ for\ } 1\le i \le n-1,
\end{array}
\]
and they are known to be counted by the Catalan number $C_{n}$ (essentially 
interpretation ($w^{4})$ in \cite{ec2}).

\lr specifications provide a known way to show the 
Wilf-equivalence (same 
counting sequence) of 3-2-1-avoiding and 3-1-2-avoiding permutations, which we review because we'll 
need it. For a valid \lr specification for $[n]$, define its \emph{minimal} 
permutation 
 \[
\begin{array}{clcccccc}
 \textrm{position} &   \,1 \ \ \ 2   & \ldots &  p_{2} & \ldots &   p_{3}  & \ldots      \\
  \textrm{entry} &  m_{1}\, \square & \ldots & \square \, m_{2}\, \square & \ldots & \square \, m_{3} \,
    \square & \ldots  
\end{array}
\]
by filling in the boxes left to right with the smallest element of 
$[n]$ not yet used. Similarly, define its \emph{maximal} 
permutation  by filling in the boxes left to right with the largest 
unused element of 
$[n]$ that does not create a new \lr entry. 
 For example, the minimal permutation 
for  $(P,M)=((1,3,6),(3,5,7))$ is $\red{3}\,1\,\red{5}\,2\,4\,\red{7}\,6$ 
and its maximal permutation is $\red{3}\,2\,\red{5}\,4\,1\,\red{7}\,6$.
Then we have the routinely-verified
characterizations
\begin{itemize}
    \item  A permutation on $[n]$ is 3-2-1-avoiding $\Leftrightarrow$ 
    it is the minimal permutation for its \lr specification

    \item  A permutation on $[n]$ is 3-1-2-avoiding $\Leftrightarrow$ 
    it is the maximal permutation for its \lr specification,
\end{itemize}
yielding a classic bijection \cite{simion-schmidt} (see also 
\cite{132diagram}) from 3-2-1-avoiding to 3-1-2-avoiding permutations 
on $[n]$.

Decomposing a permutation $\pi$ as $m_{1}L_{1}m_{2}L_{2}\ldots m_{r}L_{r}$ with
$m_{i}$ the \lr entries of $\pi$, the following characterizations are 
now
straightforward, if a little tedious, to verify.
\begin{itemize}
    \item
     A permutation $\pi$ on $[n]$ is \s $\Leftrightarrow$ (i) each 
    $L_{i}$ is \s, and (ii) the permutation obtained from $\pi$ by 
    sorting each $L_{i}$ in increasing order is 3-2-1-avoiding (i.e. 
    minimal).
   \item
    A permutation $\pi$ on $[n]$ is \m  $\Leftrightarrow$ 
    $($i\,$)$ each 
    $L_{i}$ is \m, and $($ii$\,)$ the permutation obtained from $\pi$ by 
    sorting each $L_{i}$ in decreasing order is $3$-$1$-$2$-avoiding 
    (i.e. maximal).
\end{itemize}
These characterizations, together with the above 3-2-1 to 3-1-2 bijection, 
permit a recursively defined bijection that preserves the \lr specification,
from \s to \mbox{31-4-2}-avoiding
permutations on $[n]$. The Theorem follows.

\vspace*{10mm}

\centerline{ {\large \textbf{Acknowledgment}  }}

I thank Richard Brualdi for posing a question that prompted this paper.

\vspace*{10mm}


\end{document}